\documentclass[12pt]{article}
\usepackage{}
\usepackage{mathrsfs}

\usepackage{amsmath, epsfig, cite}
\usepackage{amssymb}
\usepackage{amsfonts}
\usepackage{latexsym}
\usepackage{graphicx,ifpdf}

\newtheorem{thm}{Theorem}[section]

\newtheorem{lem}[thm]{Lemma}

\newcommand{\qed}{{\hfill\rule{4pt}{7pt}}}
\def\pf{\noindent {\it Proof.} }

\numberwithin{equation}{section}

\makeatletter \@addtoreset{equation}{section} \makeatother

\begin{document}
\rule{0cm}{1cm}
\begin{center}
{\Large\bf On reformulated zagreb indices with respect to tricyclic graphs \footnote{The first author is supported by  NNSFC (Nos. 11326216 and 11301306)
 .}
}
\end{center}
 \vskip 2mm \centerline{ Shengjin Ji$^{1}$, Xia Li$^{2}$, Yongke Qu$^
 3$
 \footnote{
E-mail addresses: jishengjin2013@163.com(S.Ji),summer08lixia@163.com(X.Li),yongke1239@163.com(Y. Qu)}}

\begin{center}
$^{1\,2}$ School of Science, Shangdong University of Technology,
\\ Zibo, Shandong 255049, China\\
$^3$ Department of Mathematics, Luoyang Normal University,\\
 Luoyang, Henan 471022, P. R. China

\end{center}

\begin{center}
{\bf Abstract}
\end{center}

{\small The authors Mili$\breve{c}$evi$\acute{c}$ et al. introduced the reformulated Zagreb indices \cite{M&N&T2004}, which is
a generalization of classical Zagreb indices of chemical graph theory.
In the paper, we characterize the extremal properties of the first reformulated Zagreb index. We first introduce some graph operations which increase or decrease
this index. Furthermore, we will determine the extremal
tricyclic graphs with minimum and maximum the first Zagreb index by these graph operations.

\vskip 3mm

\noindent {Keywords: The reformulated Zagreb index, Zagreb indices,
Tricyclic graph, edge degree, Graph operation} \vskip 3mm

\vskip 3mm \noindent {\bf AMS Classification:} 05C50, 05C35

\section{Introduction}

 Topological indices are major invariants to characterize some properties of the graph of a molecule.
One of the most important
topological indices is the well-known Zagreb indices, as a pair of molecular descriptors, introduced in \cite{T&C2000,Gutman1972}.
For a simple graph $G$, the first and second Zagreb
indices, $M_1$ and $M_2$, respectively, are defined as:

\begin{equation*}
M_1(G)=\sum_{v\in V}deg(v)^2,\ \ \
M_2(G)=\sum_{uv\in E}deg(u)\cdot deg(v).
\end{equation*}

Zagreb indices, as a pair of molecular descriptors,
 first appeared in the topological formula for the total дл- energy
of conjugated molecules that has been derived in 1972 \cite{Gutman1972}. Soon after these indices
have been used as branching indices \cite{Gutman1975}. Later the Zagreb indices found applications
in QSPR and QSAR studies\cite{T&C2000,B1997,D1999}.

Since an edge of graph $G$ transform to a corresponding vertex of the line graph $L(G)$. Motivated by the connection,
Mili$\breve{c}$evi$\acute{c}$ et al. \cite{M&N&T2004} in 2004 reformulated the Zagreb indices in terms of edge-degrees instead of vertex-degrees as:
\begin{equation*}
EM_1(G)=\sum_{e\in E}deg(e)^2,\ \ \
EM_2(G)=\sum_{e\sim f}deg(e)\cdot deg(f),
\end{equation*}where $deg(e)$ denotes the degree of the edge $e$ in $G$, which is defined as $deg(e)=deg(u)+deg(v)-2$ with $e=uv$, and $e\sim f$ means that
the edges $e$ and $f$ are adjacent, i.e., they share a common end-vertex in
$G$.

In order to exhibit our results, we introduce some graph-theoretical notations and terminology. For other undefined ones,
see the book \cite{graphBondy2008}.

Let $S_n$, $P_n$ and $C_n$ be the star, path and cycle on n vertices, respectively. Let $G=(V; E)$ be a simple undirected graph. For $v\in V(G)$ and $e\in E(G)$, let $N_G(v)$ (or N(v) for short) be the set of all neighbors of $v$
in $G$, $G-v$ be a subgraph of $G$ by deleting vertex $v$, and $G-e$ be a subgraph of $G$ by deleting edge $e$. Let $G_0$ be a nontrivial graph and $u$ be its vertex. If $G$ is obtained by $G_0$ fusing a tree $T$ at $u$. Then we say that $T$ is a \emph{subtree} of $G$ and $u$ is its \emph{root}. Let $u\circ v$ denote the fusing two vertices $u$ and $v$ of $G$. Let $S_n^m$ denote the graph obtained by connecting one pendent to $m-n+1$ others pendents of $S_n$. In addition,
we replace the sign ``if and only if'' by ``if{f}'' for short.

Recently, the upper and lower bounds on $EM_1(G)$ and $EM_2(G)$ were presented in \cite{Zhou&T2010,I&zhou2012,De2012}; Su et al. \cite{Su&Xiong2011} characterize the extremal graph properties on $EM_1(G)$ with respect to given connectivity. As some examples, we now introduce the extremal of $EM_1(G)$ among acyclic,  unicyclic, bicyclic graphs, respectively.
\begin{thm}\label{acyclic}
Let $G$ be a acyclic connected graph with order $n$. Then
\begin{equation*}
EM_1(P_n)\leq EM_1(G)\leq EM_1(S_n),
\end{equation*}while the lower bound is attached i{ff} $G\cong P_n$ and
the upper bound is attached i{ff} $G\cong S_n$.
\end{thm}

Ili$\acute{c}$ and Zhou \cite{I&zhou2012} obtained the next conclusion. Ji and Li\cite{ji&li2014} provided a shorter proof by utilizing some graph operations.

\begin{thm}\label{unicyclic}
Let $G$ be a unicyclic graph with $n$ vertices. Then
\begin{equation*}
EM_1(C_n)\leq EM_1(G)\leq EM_1(S^{n}_n),
\end{equation*}while the lower bound is attached i{ff} $G\cong C_n$ and
the upper bound is attached i{ff} $G\cong S^n_n$.
\end{thm}

In \cite{ji&li2014}, the authors also got the bound of $EM_1$ among bicyclic graphs and completely characterized the extremal graphs correspondingly.

\begin{thm}\label{bicyclic}
Let $G$ be a bicyclic graph with $n$ vertices. Then
\begin{equation*}
4n+34\leq EM_1(G)\leq n^3-5n^2+16n+4,
\end{equation*}where the lower bound is attached i{ff} $G\in\{P_n^{k,
\ell,m}:\ell\geq 3 \}\cup\{C_n(r,\ell,t):\ell\geq 3\}$ and
the upper bound is attached i{ff} $G\cong S^{n+1}_n$.
\end{thm} The latest related results on $EM_1$ refer to \cite{d&x2013,f&s2014,g&f2014,xu&das2014,zhong&xu2014}.

In this paper we characterize the extremal properties of the first reformulated Zagreb index. In Section 2 we present some graph operations which increase or decrease
$EM_1$. In Section 3, we determine the extremal
 tricyclic graphs with minimum and maximum the first Zagreb index.

\section{Some graph operations}
In the section we will introduce some graph operations, which increase or decrease
the first reformulated Zagreb index. In fact, these graph operations will play an key role in determining the extremal graphs of the first reformulated Zagreb index among all tricyclic graphs.

Now we introduce two graph operations \cite{ji&li2014} which strictly increases the first reformulated Zagreb index of a graph.

{\bf Operation I.} As shown in Fig. 1, let $uv$ be an edge of connected graph $G$ with $d_G(v)\geq 2.$ Suppose that $\{v,w_1,w_2,\cdots w_t\}$ are all the neighbors of vertex $u$ while $w_1,w_2,\cdots w_t$ are pendent vertices. If $G'=G-\{uw_1,uw_2,\cdots uw_t\}+\{vw_1,vw_2,\cdots vw_t\}$, we say that $G' $ is obtained from $G$ by \emph{Operation I}.

\vspace{10mm}

\setlength{\unitlength}{1mm}
\begin{picture}(40,20)

\put(40,20){\circle{20}}
\put(47,20){\circle*{1.3}}
\put(52,20){\circle*{1.3}}
\put(57,26){\circle*{1.3}}
\put(57,22){\circle*{1.3}}\put(57,14){\circle*{1.3}}
\put(47,20){\line(1,0){5}}\put(52,20){\line(4,5){5}}
\put(52,20){\line(2,1){5}}\put(52,20){\line(4,-5){5}}
\multiput(57,17)( 0,1.5){3}{\circle*{0.7}}

\put(45,20){\makebox(0,0){$v$ }}
\put(52,22){\makebox(0,0){$u$ }}
\put(61,26){\makebox(0,0){$w_1$ }}\put(61,22){\makebox(0,0){$w_2$ }}
\put(61,14){\makebox(0,0){$w_t$ }}
\put(40,20){\makebox(0,0){$G_0$ }}


\put(95,20){\circle{20}}
\put(102,20){\circle*{1.3}}
\put(107,20){\circle*{1.3}}
\put(112,26){\circle*{1.3}}
\put(112,22.5){\circle*{1.3}}\put(112,14){\circle*{1.3}}

\put(102,20){\line(1,0){5}}
\put(102,20){\line(4,1){10}}\put(102,20){\line(5,3){10}}
\put(102,20){\line(5,-3){10}}
\multiput(112,17)(0, 1.5){3}{\circle*{0.7}}
\put(100,20){\makebox(0,0){$v$ }}
\put(107,18){\makebox(0,0){$u$ }}
\put(116,26){\makebox(0,0){$w_1$ }}\put(116,22){\makebox(0,0){$w_2$ }}
\put(116,14){\makebox(0,0){$w_t$ }}

\put(95,20){\makebox(0,0){$G_0$}}

\put(50,9){\makebox(0,0){$G$
}} \put(94,9){\makebox(0,0){$G'$ }}
\put(75,17){\makebox(0,0){$\longrightarrow$
}}\put(75,20){\makebox(0,0){\emph{Operation I}
}}
\put(65,0){\makebox(0,0){Fig. 1 \ Graphs $G$ and $G'$ in Operation I. }}

\end{picture}
\vspace{0mm}

{\bf Operation I{I}.} As shown in Fig. 2, Let $G$ be nontrivial connected graph $G$ and $u$ and $v$ be two vertices of $G$. Let $P_\ell=v_1(=u)v_2\cdots v_\ell(=v)$ is a nontrivial $\ell$- length path of $G$ connecting vertices $u$ and $v$. If $G'=G-\{v_1v_2,v_2v_3,\cdots,v_{\ell-1}v_{\ell}\}+\{w(=u\circ v)v_1,wv_2,\cdots,wv_{\ell}\}$, we say that $G' $ is obtained from $G$ by \emph{Operation I{I}}.

\vspace{15mm}

\setlength{\unitlength}{1mm}
\begin{picture}(40,20)

\put(22,20){\circle{20}}\put(50,20){\circle{20}}
\put(32,22){\circle*{1.2}}\put(39,22){\circle*{1.2}}
\put(28.7,22){\circle*{1.2}}\put(43.3,22){\circle*{1.2}}
\put(28.7,22){\line(1,0){4}}\put(39,22){\line(1,0){5}}
\multiput(34,22)(1.5,0){3}{\circle*{0.7}}
\put(28.5,20){\makebox(0,0){$u$ }}\put(44.5,20){\makebox(0,0){$v$ }}
\put(22,20){\makebox(0,0){$H_{1}$ }}\put(50,20){\makebox(0,0){$H_2$ }}
\put(36,25){\makebox(0,0){$\overbrace{\rule{14mm}{0mm}}^{\makebox(0,2){$P_\ell$ }}$}}

\put(90,20){\circle{20}}\put(104,20){\circle{20}}

\put(97,22){\circle*{1.2}}
\put(93,30){\circle*{1.2}}\put(101,30){\circle*{1.2}}
\put(97,22){\line(-1,2){4}}\put(97,22){\line(1,2){4}}

\multiput(95.5,30)(1.5,0){3}{\circle*{0.7}}

\put(91,31){\makebox(0,0){$v_1$ }}\put(104,31){\makebox(0,0){$v_{\ell-1}$ }}
\put(97,20){\makebox(0,0){$w$ }}

\put(90,20){\makebox(0,0){$H_{1}$ }}\put(104,20){\makebox(0,0){$H_2$ }}
\put(68,20){\makebox(0,0){$\longrightarrow$ }}
\put(69,23){\makebox(0,0){\emph{Operation {I{I}}} }}
\put(45,10){\makebox(0,0){$G$ }}
\put(97,10){\makebox(0,0){$G'$ }}
\put(60,3){\makebox(0,0){Fig. 2  The graphs $G$ and $G'$
 in Operation {I{I}}. }}
\end{picture}
\vspace{-2mm}

In fact, those inverse operation of Operation I and Operation I{I} decrease $EM_1$ of a graph. According the above two graph operations, it is immediate
to get the following two results \cite{ji&li2014}.
\begin{lem}\label{operationI}
 If $G'$ is obtained from $G$ by Operation I as shown in Fig. 1, then
$$ EM_1(G)<EM_1(G').$$
\end{lem}
\begin{lem}\label{operationII}
If $G'$ is obtained from $G$ by Operation I{I} as shown in Fig. 2, then
$$EM_1(G)<EM_1(G').$$
\end{lem}

{\bf Operation I{I}{I}.} As shown in Fig. 3, let $G_0$ be a nontrivial subgraph(acyclic) of $G$ with $|G_0|=t$ which is
attached at $u_1$ in graph $G$, $x$ and $y$ be two neighbors of $u_1$ different from in $G_0$. If $G'=G-(G_0-u_1)+u_1v_2+v_2v_3+\cdots+v_{t}y$, we say that $G' $ is obtained from $G$ by \emph{Operation I{I}{I}}.

\vspace{10mm}

\setlength{\unitlength}{1mm}
\begin{picture}(40,20)

\put(25,19){\circle*{1.2}}\put(25,25){\circle*{1.2}}
\put(25,15){\circle*{1.2}}
\put(22,18){\circle*{1.2}}
\put(22,12){\circle*{1.2}}
\put(25,15){\circle*{1.2}}\put(34,15){\circle*{1.2}}
\put(37,18){\circle*{1.2}}\put(37,12){\circle*{1.2}}
\put(25,15){\line(-1,1){4}}\put(25,15){\line(-1,-1){4}}
\put(34,15){\line(1,1){4}}\put(34,15){\line(1,-1){4}}

\put(25,15){\line(1,0){9}}\put(25,15){\line(1,0){7}}
\put(25,15){\line(0,1){4}}
\multiput(22,14)( 0,1.3){3}{\circle*{0.5}}
\multiput(37,14)( 0,1.3){3}{\circle*{0.5}}
\multiput(25,21)( 0,1.3){3}{\circle*{0.5}}
\put(28,13){\makebox(0,0){$u_1$ }}
\put(28,19){\makebox(0,0){$u_2$ }}\put(28,25){\makebox(0,0){$u_t$ }}
\put(33,13){\makebox(0,0){$y$ }}

\put(50,19){\circle{8}}
\put(50,15){\circle*{1.2}}
\put(47,18){\circle*{1.2}}
\put(47,12){\circle*{1.2}}
\put(50,15){\circle*{1.2}}\put(59,15){\circle*{1.2}}
\put(62,18){\circle*{1.2}}\put(62,12){\circle*{1.2}}
\put(50,15){\line(-1,1){4}}\put(50,15){\line(-1,-1){4}}
\put(59,15){\line(1,1){4}}\put(59,15){\line(1,-1){4}}

\put(50,15){\line(1,0){9}}\put(50,15){\line(1,0){7}}
\put(50,20){\makebox(0,0){$G_0$}}

\multiput(47,14)( 0,1.3){3}{\circle*{0.5}}
\multiput(62,14)( 0,1.3){3}{\circle*{0.5}}
\put(52,13){\makebox(0,0){$u_1$ }}
\put(58,13){\makebox(0,0){$v$ }}

\put(97,18){\circle*{1.2}}
\put(97,12){\circle*{1.2}}
\put(100,15){\circle*{1.2}}\put(117,15){\circle*{1.2}}
\put(120,18){\circle*{1.2}}\put(120,12){\circle*{1.2}}
\put(105,15){\circle*{1.2}}\put(112,15){\circle*{1.2}}
\put(100,15){\line(-1,1){4}}\put(100,15){\line(-1,-1){4}}
\put(100,15){\line(1,0){5}}
\put(112,15){\line(1,0){5}}
\put(117,15){\line(1,1){4}}\put(117,15){\line(1,-1){4}}
\multiput(107,15)(1.5,0 ){3}{\circle*{0.6}}
\multiput(97,14)(0, 1.2){3}{\circle*{0.5}}
\multiput(120,14)(0, 1.2){3}{\circle*{0.5}}

\put(101,13){\makebox(0,0){$u_1$ }}
\put(105,13){\makebox(0,0){$u_2$ }}\put(112,13){\makebox(0,0){$u_t$ }}
\put(116,13){\makebox(0,0){$v$ }}

\put(30,7){\makebox(0,0){$G_1$}}

\put(55,7){\makebox(0,0){$G$
}} \put(105,7){\makebox(0,0){$G'$ }}
\put(80,15){\makebox(0,0){$\longrightarrow$
}}\put(79,18){\makebox(0,0){\emph{Operation I{I}{I}}
}}
\put(65,0){\makebox(0,0){Fig. 3 \ Graphs $G$, $G'$, $G_1$ in
Operation I{I}{I}. }}

\end{picture}
\vspace{3mm}

\begin{lem}\label{operationIII}
Let $G$ and $G'$ be two graphs as shown in Fig. 3. Then we have
$$ EM_1(G)>EM_1(G')$$
\end{lem}
\pf According to the inverse of Operation I, as shown in Fig. 3, there is a graph $G_1$ such that
$EM_1(G)\geq EM_1(G_1)$. In order to show the conclusion, we now just to verify the following Inequality:
\begin{equation}\label{ivcomparation}
EM_1(G_1)>EM_1(G').
\end{equation} By means of the definition of $EM_1$, we have
\begin{equation*}
\begin{split}
EM_1(G_1)-EM_1(G')&>d^2_{G_1}(u_{1}u_2)
                    +d^2_{G_1}(u_1v)+d^2_{G_1}(u_{t-1}u_t)\\
                   &\quad-[ d^2_{G'}(u_{1}u_2)+d^2_{G'}(u_{t-1}u_t)
                    +d^2_{G'}(u_tv)]\\
                    &=d^2_{G_1}(u_1)+(d_{G_1}(u_1)+d_{G_1}(v))^2-
                    (d_{G_1}(u_1)-1)^2-d^2_{G_1}(v)-1.\\
                    &=d^2_{G_1}(u_1)+2d_{G_1}(u_1)d_{G_1}(v)+2d_{G_1}(u_1)-2>0
\end{split}
\end{equation*} Therefore, the Ineq. (\ref{ivcomparation}) holds. Then we
finish the proof. \qed

{\bf Operation I{V} .} Let $G_0$ be a nontrivial connected graph and
 $u$ and $v$ are two vertices in $G_0$ with
 $d_{G_0}(u)=x, d_{G_0}(v)=y$ and $N_{G_0}(u)\supseteq N_{G_0}(v)$. Let $G$ be the graph obtained by attaching $S_{k+1}$
and $S_{\ell+1}$ at the vertices $u$ and $v$ of $G_0$, respectively. If
 $G'$ is the graph obtained by delating the $\ell$ pendent vertices at $v$ in $G$ and connecting them to $u$ of $G$, depicted in Fig. 4,
We say that $G'$ is obtained from
$G$ by \emph{Operation I{V}}.

\vspace{15mm}

\setlength{\unitlength}{1mm}
\begin{picture}(40,20)

\put(45,20){\circle{20}}
\put(51,24){\circle*{1.2}}
\put(51,16){\circle*{1.2}}
\put(56,28){\circle*{1.2}}
\put(56,22){\circle*{1.2}}\put(56,12){\circle*{1.2}}
\put(56,18){\circle*{1.2}}
\put(51,24){\line(5,4){5}}\put(51,24){\line(3,-1){5}}
\put(51,16){\line(5,-4){5}}\put(51,16){\line(3,1){5}}
\multiput(56,23.5)( 0,1.5){3}{\circle*{0.5}}
\multiput(56,13.5)( 0,1.5){3}{\circle*{0.5}}
\put(49,24){\makebox(0,0){$u$ }}\put(49,16){\makebox(0,0){$v$ }}
\put(59,28){\makebox(0,0){$u_1$ }}\put(59,22){\makebox(0,0){$u_k$ }}
\put(59,18){\makebox(0,0){$v_1$ }}\put(59,12){\makebox(0,0){$v_\ell$ }}
\put(45,20){\makebox(0,0){$G_0$}}

\put(95,20){\circle{20}}
\put(101,24){\circle*{1.2}}
\put(107,23){\circle*{1.2}}
\put(107,25){\circle*{1.2}}
\put(107,30){\circle*{1.2}}\put(107,18){\circle*{1.2}}

\put(101,24){\line(5,1){6}}\put(101,24){\line(5,-1){6}}
\put(101,24){\line(1,1){6}}\put(101,24){\line(1,-1){6}}
\multiput(107,26)(0, 1.2){3}{\circle*{0.5}}
\multiput(107,19)(0, 1.2){3}{\circle*{0.5}}

\put(99,24){\makebox(0,0){$u$ }}
\put(110,30){\makebox(0,0){$u_1$ }}\put(110,25.5){\makebox(0,0){$u_k$ }}
\put(110,22.5){\makebox(0,0){$v_1$ }}\put(110,18){\makebox(0,0){$v_\ell$ }}

\put(95,20){\makebox(0,0){$G_0$}}

\put(50,7){\makebox(0,0){$G$
}} \put(90,7){\makebox(0,0){$G'$ }}
\put(75,17){\makebox(0,0){$\longrightarrow$
}}\put(75,20){\makebox(0,0){\emph{Operation I{V}}
}}
\put(65,0){\makebox(0,0){Fig. 4 \ $G$ and $G'$ in Operation I{V}. }}
\end{picture}
\vspace{-3mm}

\begin{lem}\label{operationiv}
If $G'$ is obtained from $G$ by Operation V as shown in Figure. 4. Then
$$EM_1(G)<EM_1(G').$$
\end{lem}
\pf Note that $d_{G_0}(u)=x,$ and $d_{G_0}(v)=y>0$, meanwhile $N_{G_0}(u)\supseteq N_{G_0}(v)$. By the definition of $EM_1$, we have
 \begin{equation*}
\begin{split}
EM_1(G')-EM_1(G)&>\sum^{k}_{i=1}[d^2_{G'}(uu_i)-d^2_G(uu_i)]
+\sum^{\ell}_{i=1}[d^2_{G'}(uv_i)-d^2_G(vv_i)]\\
               &\ +\sum_{w\in N_{G_0}(v)}[d_{G'}^2(uw)+d_{G'}^2(vw)]
               -\sum_{w\in N_{G_0}(v)}[d_{G}^2(uw)+d_{G}^2(vw)]\\
                &=(k+\ell)(k+\ell+x-2)^2-k(k+x-2)^2-\ell(\ell+y-2)^2\\
                &\,+\sum_{w\in N_{G_0}(v)}[(k+\ell+x-d_{G_0}(w)-2)^2+
                (y+d_{G_0}(w)-2)^2]\\&\,-\sum_{w\in N_{G_0}(v)}[(k+x+d_{G_0}(w)-2)^2+(\ell+y+d_{G_0}(w)-2)^2]\\
                &>2\ell(x+k-y)>0.
\end{split}
\end{equation*}

 So the result follows. \qed

 As the above exhibited, Operation I{I}{I} strictly decrease
 the $EM_1$ of a graph; while all of Operation I, Operation I{I} and Operation I{V} strictly increase the $EM_1$ of a graph.

\section{Main results}
In the section, we will characterize the extremal graph with respect to $EM_1$ among all tricyclic graphs by some graph operations.

\begin{figure}[ht]
\begin{center}
\includegraphics[width=10cm]{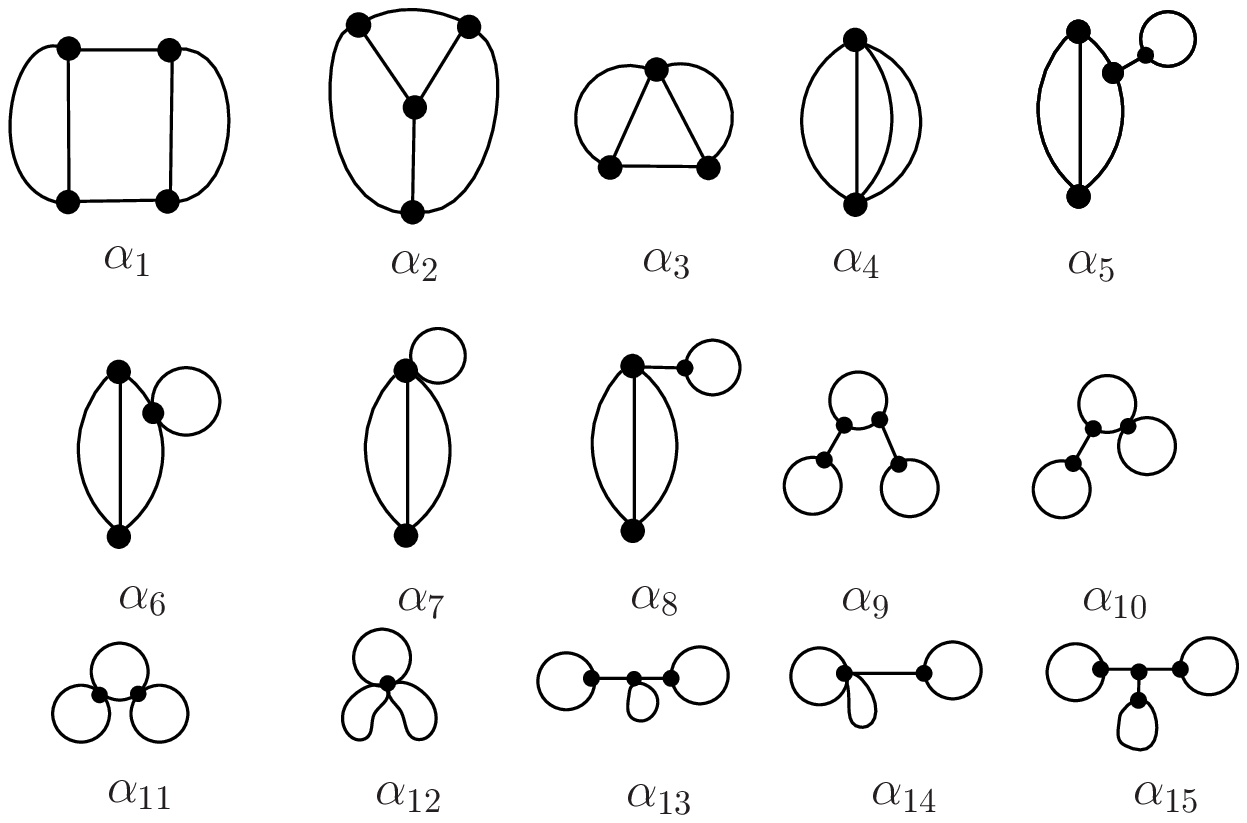}
\end{center}
\end{figure}
\vspace{-13mm}
\begin{center}
Fig. 5 The graphs in $\mathscr{C}_n^0$
\end{center}
\vspace{-5mm}

For convenience, we now define some notations which will be using in the sequel.
Denote by $\mathscr{C}_n$ the set of all connected tricyclic graphs with order $n$. For any tricyclic graph $G$, the subgraph which is obtained by deleting all pendents of $G$ is referred as a \emph{brace} of $G$.
Let $\mathscr{C}_n^0$ be the set of all braces of tricyclic graphs as pictured in Fig. 5, and $\mathscr{C}_n^1$ denote the set of these tricyclic graphs  shown in Fig. 6. Moreover, $\mathscr{C}_n^2$ denote the set of these tricyclic graphs depicted in Fig. 7.

\begin{figure}[ht]
\begin{center}
\includegraphics[width=10cm]{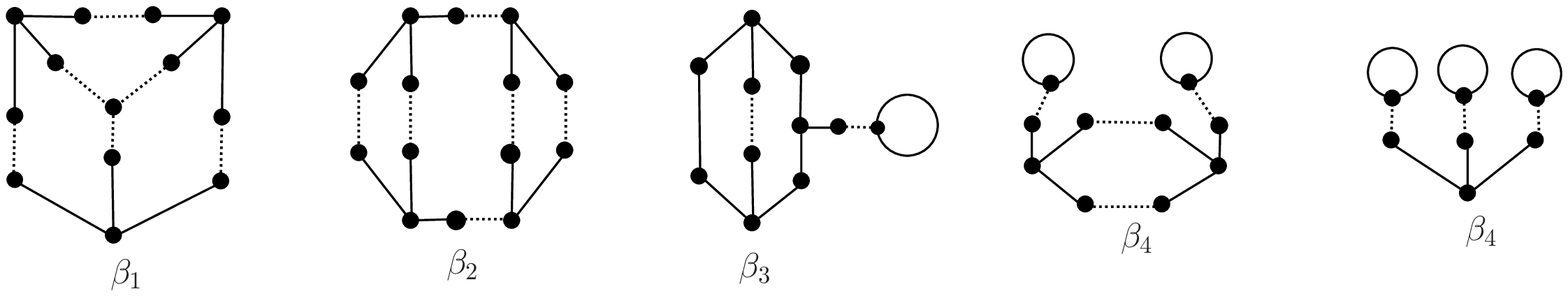}
\end{center}
\end{figure}
\vspace{-11mm}
\begin{center}
Fig. 6 The graphs in $\mathscr{C}_n^1$
\end{center}

\vspace{2mm}

\setlength{\unitlength}{1mm}
\begin{picture}(40,20)
\put(20,10){\circle*{1.2}}\put(23,10){\circle*{1.2}}
\put(26,10){\circle*{1.2}}\put(29,10){\circle*{1.2}}
\put(24.5,17){\circle*{1.2}}
\put(24.5,21){\circle*{1.2}}\put(22.5,21){\circle*{1.2}}
\put(28.5,21){\circle*{1.2}}
\put(24.5,17){\line(-2,-3){5}}\put(24.5,17){\line(2,-3){5}}
\put(24.5,17){\line(-1,-4){1.8}}\put(24.5,17){\line(1,-4){1.8}}
\put(20,10){\line(1,0){10}}
\put(24.5,17){\line(0,1){4}}
\put(24.5,17){\line(-1,2){2}}\put(24.5,17){\line(1,1){4}}
\multiput(25.5,21)(1.0,0){3}{\circle*{0.5}}



\put(35,10){\circle*{1.2}}\put(38,10){\circle*{1.2}}
\put(41,10){\circle*{1.2}}
\put(43,10){\circle*{1.2}}
\put(46,10){\circle*{1.2}}\put(40,17){\circle*{1.2}}

\put(40,21){\circle*{1.2}}\put(38,21){\circle*{1.2}}
\put(44,21){\circle*{1.2}}

\put(40,17){\line(-2,-3){5}}
\put(40,17){\line(-1,-3){2.5}}\put(40,17){\line(2,-5){3}}
\put(40,17){\line(1,-6){1.2}}\put(40,17){\line(5,-6){6}}
\put(35,10){\line(1,0){6}}\put(43,10){\line(1,0){3}}

\put(40,17){\line(0,1){4}}
\put(40,17){\line(-1,2){2}}\put(40,17){\line(1,1){4}}
\multiput(41,21)(1.0,0){3}{\circle*{0.5}}


\put(52,10){\circle*{1.2}}\put(54.5,10){\circle*{1.2}}
\put(57,10){\circle*{1.2}}\put(61,10){\circle*{1.2}}
\put(63.5,10){\circle*{1.2}}\put(66,10){\circle*{1.2}}

\put(59,17){\circle*{1.2}}
\put(59,21){\circle*{1.2}}\put(57,21){\circle*{1.2}}
\put(63,21){\circle*{1.2}}

\put(59,17){\line(-1,-3){2.5}}\put(59,17){\line(-2,-3){5}}
\put(59,17){\line(-1,-1){7.5}}
\put(59,17){\line(1,-3){2.5}}\put(59,17){\line(2,-3){5}}
\put(59,17){\line(1,-1){7.5}}

\put(52,10){\line(1,0){2}}\put(57,10){\line(1,0){4}}
\put(63.5,10){\line(1,0){2}}

\put(59,17){\line(0,1){4}}
\put(59,17){\line(-1,2){2}}\put(59,17){\line(1,1){4}}
\multiput(60,21)(1.0,0){3}{\circle*{0.5}}


\put(75,10){\circle*{1.2}}\put(79,10){\circle*{1.2}}
\put(81,10){\circle*{1.2}}\put(85,10){\circle*{1.2}}
\put(78,14){\circle*{1.2}}\put(82,14){\circle*{1.2}}
\put(80,18){\circle*{1.2}}

\put(80,22){\circle*{1.2}}\put(78,22){\circle*{1.2}}
\put(84,22){\circle*{1.2}}
\put(80,18){\line(-1,-2){2}}
\put(80,18){\line(1,-2){2}}
\put(78,14){\line(-4,-5){3.5}}\put(78,14){\line(1,0){4}}
\put(78,14){\line(1,-4){1.2}}
\put(82,14){\line(4,-5){3.5}}
\put(82,14){\line(-1,-4){1.2}}\put(75,10){\line(1,0){4}}
\put(81,10){\line(1,0){4}}

\put(80,18){\line(0,1){4}}
\put(80,18){\line(-1,2){2}}\put(80,18){\line(1,1){4}}
\multiput(81,22)(1.0,0){3}{\circle*{0.5}}

\put(95,14){\circle*{1.2}}
\put(99,10){\circle*{1.2}}\put(101,14){\circle*{1.2}}
\put(99,18){\circle*{1.2}}\put(103,14){\circle*{1.2}}

\put(99,22){\circle*{1.2}}\put(97,22){\circle*{1.2}}
\put(103,22){\circle*{1.2}}

\put(99,10){\line(-1,1){4}}\put(99,10){\line(1,1){4}}
\put(99,18){\line(-1,-1){4}}\put(99,18){\line(1,-1){4}}
\put(99,18){\line(0,-1){8}}
\put(101,14){\line(-1,2){4}}\put(101,14){\line(-1,-2){2}}

\put(99,18){\line(0,1){4}}
\put(99,18){\line(-1,2){2}}\put(99,18){\line(1,1){4}}
\multiput(100,22)(1.0,0){3}{\circle*{0.5}}



\put(111,14){\circle*{1.2}}
\put(115,10){\circle*{1.2}}
\put(115,18){\circle*{1.2}}\put(119,14){\circle*{1.2}}

\put(115,22){\circle*{1.2}}\put(113,22){\circle*{1.2}}
\put(119,22){\circle*{1.2}}

\put(115,10){\line(-1,1){4}}\put(115,10){\line(1,1){4}}
\put(115,18){\line(-1,-1){4}}\put(115,18){\line(1,-1){4}}
\put(115,18){\line(0,-1){8}}\put(111,14){\line(1,0){8}}

\put(115,18){\line(0,1){4}}
\put(115,18){\line(-1,2){2}}\put(115,18){\line(1,1){4}}
\multiput(116,22)(1.0,0){3}{\circle*{0.5}}



\put(25,7){\makebox(0,0){$\gamma_1$ }}
\put(40,7){\makebox(0,0){$\gamma_2$ }}
\put(60,7){\makebox(0,0){$\gamma_3$}}
\put(80,7){\makebox(0,0){$\gamma_4$ }}
\put(100,7){\makebox(0,0){$S_n^{n+2}(\gamma_5)$ }}
\put(118,7){\makebox(0,0){$S_n^{K_4}(\gamma_6)$ }}
\put(65,0){\makebox(0,0){Fig. 7 \ Some graphs using in the later proof. }}

\end{picture}
\vspace{3mm}

We next introduce the extremal graphs with respect to $EM_1$ on tricyclic graphs.

\begin{thm}\label{minitri}
Let $G$ be a tricyclic graph with order $n$. Then
\begin{equation*}
4n+68\leq EM_1(G),
\end{equation*}where the equality holds i{ff} $G\in \mathscr{C}_n^1$.
\end{thm}
\pf Let $G$ be a connected tricyclic graph. By Lemma \ref{operationII}, $G$ can be converted to the one of the fifteen graphs without any pendent as shown in Fig. 5. In other words, there exists a graph $\alpha_i\in \mathscr{C}_n^0$ (\,$i\leq 15$\,) such that $EM_1(G)\geq EM_1(\alpha_i)$ in terms of Lemma \ref{operationIII}, for any given graph $G$. By directed calculation, we obtain that
\begin{equation*}
EM_1(\beta_i)=4n+68, \ \text{for}\  i=1,2,\cdots,5.
\end{equation*}Then the proof complete. \qed

\begin{thm}\label{maxitri}
Let $G$ be a tricyclic graph with order $n$. Then
\begin{equation*}
EM_1(G)\leq n^3-5n^2+20n+32,
\end{equation*}where the equality holds i{ff} $G\cong S_n^{n+2}$ or $ S_n^{K{_4}}$.
\end{thm}
\pf For a given connected tricyclic graph $G$, with repeated Operation I{I} and Operation I{V}, it can be converted to the one of the five graphs as shown in Fig. 6. That is to say, For any tricyclic graph $G$ with order $n$, there exists $\gamma_i\in\mathscr{C}_n^2$(\,$i\leq 6$\,) such that $EM_1(G)\leq EM_1(\gamma_i)$ in views of Lemma \ref{operationII} and Lemma \ref{operationiv}.

By indirected calculating, we get
\begin{equation*}
\begin{split}
 EM_1(\gamma_1)&=n^3-5n^2+16n+18,\, EM_1(\gamma_2)=n^3-5n^2+20n-10,\\
 EM_1(\gamma_3)&=n^3-5n^2+20n+2,\,  EM_1(\gamma_4)=n^3-9n^2+32n+60,\\
EM_1(S_n^3)&=n^3-5n^2+20n+32,\, EM_1(S_n^{K_4})=n^3-5n^2+20n+32.
\end{split}
\end{equation*} Therefore we complete the proof. \qed

Together Theorem \ref{minitri} with Theorem \ref{maxitri}, the main result is shown.

\end{document}